\documentclass[12pt]{article}

\usepackage{amsfonts,amsmath, amssymb}
\usepackage{hyperref}
\usepackage{tikz}

\def\CP{\C\mathrm{P}}
\def\CH{\C\mathrm{H}}
\def\HP{\HH\mathrm{P}}
\def\OP{\OO\mathrm{P}}
\def\OH{\OO\mathrm{H}}

\def\[#1{\begin{equation}\label{#1}}
\def\]{\end{equation}}
\def\abs#1{\left|#1\right|}

\def\C{\mathbb{C}}
\def\curl{\mathop{\mathrm{curl}}}
\def\div{\mathop{\mathrm{div}}}
\def\diff#1#2{\left.\fracd{d}{d#1}\right|_{#1=#2}}
\def\diffs#1#2{\left.\frac{d^2}{d#1^2}\right|_{#1=#2}}

\def\Diffs#1#2{\left.\fracd{D^2}{d#1^2}\right|_{#1=#2}}
\def\d{\mathrm{d}}
\def\dss{\displaystyle}
\def\fracd#1#2{\frac{\raisebox{2pt}{$\dss#1$}}{\raisebox{-5pt}{$\dss#2$}}}
\def\grad{\mathrm{grad}}
\def\HH{\mathbb{H}}
\def\Hom{\mathop{\mathrm{Hom}}}
\def\id{\mathrm{id}}
\def\K{\mathbb{K}}
\def\l({\left(}
\def\N{\mathbb{N}}
\def\norm#1{\left\|#1\right\|}
\def\OO{\mathbb{O}}
\def\P{\mathcal{P}}
\def\pmx#1{\begin{pmatrix}#1\end{pmatrix}}
\def\proof{{\bf Proof: }}
\def\proofend{\hfill$\bullet$\par}
\def\qLeftrightarrow{\quad\Longleftrightarrow\quad}
\def\qmbox#1{\quad\mbox{#1}\quad}
\def\R{\mathbb{R}}
\def\r){\right)}
\def\set#1#2{\left\{#1\ \vline\hspace{110sp}\vline\hspace{110sp}\vline\ #2\right\}}
\def\sg#1{\left[#1\right]}
\def\spl#1{\begin{equation*}\begin{split}#1\end{split}\end{equation*}}
\def\spll#1#2{\begin{equation}\label{#1}\begin{split}#2\end{split}\end{equation}}
\def\spn#1{\mathop{\mathrm{span}} #1}
\def\SS{\mathrm{S}}
\def\st#1{\left\{#1\right\}}
\def\vol{\mathbf{vol}\,}

\def\KP{\K\mathrm{P}}
\def\GM{G_r(\KP^n)}
\def\VV{\overline{V}}
\def\pp{\overline{p}}
\def\upin{\rotatebox[origin=c]{90}{$\in$}}

\def\lk#1{\mathop{\mathrm{lk}}(#1)}

\def\intd{\int\displaylimits}
\def\Jac{\mathrm{Jac}_x(y,s)}

\newtheorem{thm}{Theorem}
\numberwithin{thm}{section}

\begin{document}\sloppy

\def\email#1{email: \texttt{#1}}
\def\affiliation#1{\noindent #1}

\title{Linking Integrals on Rank-1-Symmetric Spaces}

\author{Stefan Bechtluft-Sachs and Evangelia Samiou}

\maketitle
\begin{abstract}
In symmetric spaces of rank one we obtain a right inverse of the Cartan differential on exact differential forms as an integral operator. As a corollary we extend Gauss' linking integral beyond space forms.
\end{abstract}
\bigskip

\noindent{\it Keywords: Linking number, Gauss linking integral, Cartan's Formula, Rank one symmetric space}
\medskip

\noindent{\it AMS subject classifications:
58A12
; 53C35
; 53Z05
}


\def\L{\mathfrak{L}}

\section{Introduction}
The linking number of two disjoint nullhomologous closed oriented submanifolds
$K^k$, $L^l$ of complementary dimensions in an oriented manifold $M^{k+l+1}$ is the primary obstruction to disentangling them.
In Maxwell's words, it is ``the number of turns that one embraces the other in the same direction'', \cite{rini}.
The purpose of this paper is to derive integral formulas for the linking number in symmetric spaces of rank $1$, analogous to  Gauss' linking integral \eqref{gauss}.
\medskip

A consequence of the Biot-Savart Law of magnetostatics is that the linking number of two loops in Euclidean $3$-space determines the amount of work required to move a ficticious magnetic monopole along one loop in the magnetic field produced by a stationary current in the other loop.
This observation is possibly at the origin of Maxwell's rediscovery of Gauss' formula
\[{gauss}\lk{K,L}=\intd_K\intd_L\frac{\det(\d K_x\wedge(x-y)\wedge\d L_y)}{\norm{x-y}^n\vol\SS^{n-1}}\ d^kx\ d^l y\ ,\]
for the linking number in $\R^n$, 
where $\d K_x$, $\d L_y$ denote the volume elements of $K$ and $L$ at $x\in K$ and $y\in L$ respectively, \cite{rini}, \cite{epple}.
\medskip

The linking number of $K^k,L^l\subset M^{k+l+1}$ is the intersection number of $K^k$ with a $(l+1)$-chain in the surrounding manifold $M$ whose boundary is a fundamental cycle for $L^l$.
Translating this to Hodge-deRham theory gives the linking number as an integral
$$\lk{K,L}=\intd_K\d^{-1}\eta_L $$
where $\eta_L$ is a Poincar\'e dual of $L$ and $\d^{-1}\eta_L$ denotes any differential form $\alpha$ with $\d\alpha=\eta_L$, \cite{botu}.
Therefore, our main technical aim in deriving formulas like \eqref{gauss} in more general ambient spaces is the construction of a right inverse $\d^{-1}$ of the Cartan differential $\d$.
\medskip

In $\R^n$, such a right inverse is well-known and given by the Riesz-potential: If $\omega$ is a closed compactly supported $(k+1)$-form on $\R^n$, then the $k$-form $\mu$, defined at $x\in\R^n$ by 
\[{dinvrn}\mu_x(v_1,\ldots,v_k)=\intd_{\R^n}\frac{\omega_y(x-y,v_1,\ldots,v_k)}{\norm{x-y}^n\vol\SS^{n-1}}\ d^n y\ ,\]
solves both $\d\mu=\omega$ and $\d^*\mu=0$.
In $\R^3$ these two equations are equivalent to the magnetostatic part of Maxwell's equations
\spll{mstmw}{
\d\mu=\omega\qLeftrightarrow\curl\mu&=\omega\qmbox{``Ampere's Law''}\\
\d^*\mu=0\qLeftrightarrow\div\mu&=0\qmbox{``no monopoles''}
}
for the stationary magentic field $\mu$ caused by a current density $\omega$. In this context formula \eqref{dinvrn} is known as the Biot-Savart Law. In general, a Biot-Savart operator is defined in \cite{ctg} as an operator $\mathrm{BS}$ so that $\mu=\mathrm{BS}\,\omega$ solves \eqref{mstmw} for any closed differential form $\omega$.

In \cite{ctg}, \cite{dtg}, \cite{tk}, DeTurck, Gluck and Cantarella construct such an operator and derive the corresponding linking integrals in 3-dimensional space forms, see also \cite{SS} for negatively curved symmetric spaces. These results depend on the harmonicity of rank-1-symmetric spaces.
\medskip

However in order to compute the linking number from the Poincar\'e duals of the linked submanifolds, one does not need to solve the second equation $\d^*\mu=0$.
This means that the right inverses $\d^{-1}$ we will construct in this paper will give magnetic fields which may have monopoles. 
\medskip

There are some other methods to obtain linking integrals. 
In \cite{dtg3}, also see \cite{kp}, DeTurck and Gluck compute the linking number in spheres of arbitrary dimension as the degree of a map from the join of the two submanifolds.

In \cite{SV} Shonkwiler and Vela-Vick consider visible submanifolds of euclidean space, i.e. hypersurfaces intersecting once or not at all with any ray emanating from the origin. 
The linking integral in such a hypersurface can be reduced to a linking integral in the ambient euclidean space.
For round spheres, this yields the same linking integral as \cite{dtg3}, and also coincides with the formula we get in Example \ref{exlksn}.
\bigskip

The paper is organised as follows. In section \ref{sec:prel}, we describe the general form of the integral kernel of a right inverse of the Cartan differential. We also recall the main technical tool for this, the Cartan formula.
Theorem \ref{musymsp} in section \ref{sec:had} contains the formula for $\d^{-1}$ in the negatively curved symmetric spaces. Section \ref{cross} is concerned with compact symmetric spaces of rank 1. The somewhat simpler case of spheres is treated in section \ref{sec:sn}, and the projective spaces are in \ref{sec:cpn}.
Finally, in section \ref{sec:link}, we establish the connection between the integral kernel of $d^{-1}$ and linking integrals. We then compute a few examples.

\section{Preliminaries}\label{sec:prel}
\subsection{The general form of a right inverse of the Cartan differential}\label{sec:gfm}
Throughout, $M$ will be an oriented connected complete Riemannian manifold, and almost always a symmetric space of rank 1. 
Our aim is to compute a right inverse $\d^{-1}$ to the Cartan differential
$$\d\colon\Omega^{k}(M)\to\Omega^{k+1}(M)$$
in the form of an integral operator
$$\d^{-1}\colon d\Omega^{k}(M)\to\Omega^k(M)\qmbox{,}\d^{-1}\omega_x=\intd_M K^*_{yx}\omega_y\ dy\qmbox{,}\omega\in\Omega^{k+1}(M)\ ,$$
where integration is always with respect to the Riemannian metric on $M$.
In general the kernel $K^*$ of such an operator is a section in the exterior homomorphism bundle,
$$K^*\in\Gamma\Hom(\Lambda^{k+1}TM^*,\Lambda^kTM^*)$$
over $M\times M$, induced by a section
$K\in\Gamma\Hom(\Lambda^kTM,\Lambda^{k+1}TM)$. Thus
$$K^*_{yx}\omega_y=\omega_y\circ K_{yx}$$
and $K_{yx}\in\Hom(\Lambda^kT_xM,\Lambda^{k+1}T_yM)$ for $x,y\in M$.
In the cases we are interested in here, for every point $x\in M$ we will have a unique minimizing geodesic $\gamma_{yx}$ from $x$ to $y$, for almost all $y\in M$. By means of the parallel transport
$$\P_{yx}\colon T_xM\to T_yM$$
along this geodesic we can rewrite the integral kernel $K$ as
\[{ptf}K_{yx}=\Lambda^{k+1}\P_{yx}\circ k_{yx}\] 
with $k_{yx}\in\Hom(\Lambda^kT_xM,\Lambda^{k+1}T_xM)$, acting on the exterior powers at the same point $x$. 
This is the parallel-transport form of the kernel, in the terminology of \cite{ctg}, where also other possibilities of relating tangent spaces at different points are considered, for instance the left translation in a Lie group.
\medskip

The situations we will consider in this paper admit a further specialization of the kernel function.
Let $T_{yx}=\gamma'_{yx}(0)\in T_xM$ denote the unit vector at $x$ pointing to $y$ and let $\iota_{T_{yx}}$ be the contraction (see \eqref{defix}). The kernel we will construct is of the form
$$k_{yx}(v_1\wedge\cdots\wedge v_k)=\lambda(d(x,y))T_{yx}\wedge \L_{yx} v_1\wedge\cdots\wedge \L_{yx} v_k$$
with a section $\L_{yx}\in\Hom(T_xM,T_xM)$ of the endomorphisms of the tangent bundle and a function $\lambda$ depending on the distance $d(x,y)$ only.

Thus the inverse of the Cartan differential will be of the form
\[{mfdinv}(\d^{-1}\omega)_x=\intd_M \lambda(d(x,y))\L_{yx}^*\iota_{T_{yx}}\P_{yx}^*\omega_y\ dy\ ,\]
i.e.
$$(\d^{-1}\omega)_x(v_1\wedge\cdots\wedge v_k)=\intd_M \lambda(d(x,y))\omega_y(\P_{yx}(T_{yx}\wedge \L_{yx}(v_1\wedge\cdots\wedge v_k)))\ dy\ .$$

If $M$ is of constant curvature ($S^n$, $\R^n$, or $H^n$) then one can choose $\L_{yx}$ to be the identity. In these cases up to parallel transport the integral kernel of $\d^{-1}$ is a scalar function of the distance only.

\subsection{Cartan's Formula}\label{sec:cf} We will use Cartan's formula which states that, for a vector field $T$ on $M$, the Lie derivative of a form $\omega\in\Omega^{k+1}(M)$ is given by
\[{cmf}\mathcal{L}_T\omega:=\diff{t}{0}\phi_{T,t}^*\omega = \iota_T \d\omega + \d\iota_T\omega\ , \]
where
\[{defix}(\iota_T\omega)(v_1,\ldots,v_k)=\omega(T,v_1,\ldots,v_k)\]
denotes the contraction of $\omega$ with $T$ and $\phi$ is the local flow generated by $T$, see \cite{P} or \cite{B}. This means that 
$$\phi\colon U\to M\qmbox{where} M\times\st{0}\subset U\overset{\text{open}}{\subset} M\times\R$$
and $T_x=\diff{t}{0}\phi_t(x)$.
If $\omega$ is closed, i.e. $\d\omega=0$, then \eqref{cmf} simplifies to
$$\diff{t}{s}\phi_t^*\omega=\d\iota_T\phi_s^*\omega=\d\phi_s^*\iota_T\omega \ .$$ If   $U\cap \st{x}\times\R\supset\st{x}\times[0,b(x))$, for some smooth function $b\colon M\to\R^+\cup\st{\infty}$ so that $(\phi_{b(x)}^*\omega)_x=0$, then we can integrate and get $\omega=\d\mu$ with $\mu\in\Omega^k(M)$ defined by
\[{flowmain}\mu_x=-\intd_0^{b(x)}(\phi_s^*\iota_T\omega)_x\ ds\ .\]

\section{Hadamard Manifolds, Non-Compact Symmetric Spaces of Rank 1}\label{sec:had}
\subsection{The gradient flow of Busemann functions}
Let $M$ be a complete simply connected Riemannian manifold with non-positive sectional curvature.
A Busemann function $h_p$ of a point $p$ in the ideal boundary $\partial M$ is a function of the form
$$h_p(x)=\lim_{t\to\infty}(t-d(x,\gamma(t)))$$
where $\gamma$ is a geodesic representing $p$.
By \cite{bgs}, Lemma 3.4, Busemann functions $h_p$ on Hadamard manifolds are continuously differentiable and their gradient field has norm 1.
\medskip

Let $\phi_p$ be the gradient flow of a Busemann function $h_p$ of $p\in\partial M$. Then
$$T_{px}=\diff{t}{0}\phi_{p,t}x=\grad_x h_p\qmbox{,}\norm{T_p}=1 $$
and the curves
$$\gamma_{px}(t)=\phi_{p,t}(x)$$
are unit speed geodesics.
\medskip

Let $\omega\in\Omega^{k+1}(M)$ be closed of compact support. We define $\mu_p\in\Omega^k(M)$ by
\[{muhadrw}\mu_{px}:=-\intd_0^\infty\iota_{T_{px}}\phi_{p,s}^*\omega\ ds\]
Since the support of $\omega$ is compact, we have $(\phi_{p,b_x}^*\omega)_x=0$ for $b_x$ sufficiently large.
From section \ref{sec:cf} we have
$$\d\mu_p=\omega\ .$$
The differential of the flow $\phi_p$ can be expressed in terms of parallel transport and the Jacobi operator along the orbits of $\phi_p$. For $x\in M$, $v\in T_xM$, and $p\in\partial M$, 
$$\st{\phi_{p,t}(\exp_x(sv))}_{s\in\R}$$
is a geodesic variation of the geodesic $\gamma_{px}$. The variational vector field
$$d_x\phi_{p,t}v=\diff{s}{0}\phi_{p,t}(\exp_x(sv))$$
is therefore a Jacobi field along $\gamma_{px}$, invariant under the flow $\phi_p$, and bounded as $t\to\infty$.
Let $J_p\in\Gamma\Hom(TM,TM)$ denote the field of the Jacobi operators along the orbits of the flow $\phi_p$, i.e. for any $v\in TM$, we have $J_pv:=R_{T_p,v}T_p$. 
By the Jacobi equation along $\gamma_{px}$,
\spll{diffeqjac}{
\Diffs{t}{s}d_x\phi_{p,t}&=\P_{\phi_{p,s}(x),x}\diffs{t}{s}\P_{\phi_{p,t}(x),x}^{-1}d_x\phi_{p,t}\\
&= J_pd_x\phi_{p,s} \in\Hom(T_xM,T_xM)\ .
}

\subsection{Negatively Curved Symmetric Spaces}\label{secnegcurv}
In symmetric spaces the curvature tensor $R$ is parallel, hence $J_p$ commutes with the parallel transport, and \eqref{diffeqjac} becomes a ordinary differential equation for $\P_{\phi_{p,t}(x),x}^{-1}d_x\phi_{p,t}$ with constant coefficients,
\[{diffeqjacsym}\diffs{t}{s}\P_{\phi_{p,t}(x),x}^{-1}d_x\phi_{p,t} = J_p\P_{\phi_{p,s}(x),x}^{-1}d_x\phi_{p,s}\]
which now involves the Jacobi operator $J_p$ only at the point $x$.
In symmetric spaces with non-positive curvature the Jacobi operator is nonnegative. Hence the solution of \eqref{diffeqjacsym} which is bounded as $t\to\infty$ is
\[{Lnpcurv}\L_{\phi_{p,t}(x),x}:=\P_{\phi_{p,t}(x),x}^{-1} d_x\phi_{p,t}=e^{-t\sqrt{J_{px}}}\]

In the symmetric spaces with negative sectional curvature this operator has a simple eigenvalue $0$ corresponding to the eigenvector $T_{px}$. We will normalize the metric so that the other two eigenvalues will be 4 and 1.
The multiplicities of the eigenvalue $4$ are 0 in the case of hyperbolic space $H^m$, 1 in the case of $\CH^n$, 3 in the case of $\HH H^n$ and 7 in the Cayley hyperbolic plane $\OH^2$.
We can thus make \eqref{Lnpcurv} more explicit over a diagonalizing basis for $J_{px}$ to get
\spll{Lnegcurv}{\L_{\phi_{p,t}(x),x}&:=P_{\phi_{p,t}(x),x}^{-1} d_x\phi_{p,t}=e^{-t\sqrt{J_{px}}}\\
&=\pmx{
1&0    &0\\
0&e^{-2t}\id_m   &0\\
0&0&e^{-t}\id_{\dim M-m-1}
}
}
with $m=0,1,3,7$.
We can thus rewrite \eqref{muhadrw},
\[{muhadqq}\mu_{px}:=-\intd_0^\infty\L_{\phi_{p,s}(x),x}^*\iota_{T_{px}}\P_{\phi_{p,s}(x),x}^*\omega\ ds\ .\]
\medskip

For the puprose of deriving a formula invariant under the isomoetries of the ambient space, we eliminate $p$ by averaging over the ideal boundary $\partial M=\SS^{n-1}$,
$$\mu_x=\frac{-1}{\vol{\SS^{n-1}}}\intd_{\SS^{n-1}}\intd_0^\infty L_{\phi_{p,s}(x),x}^*\iota_{T_{px}}\P_{\phi_{p,s}(x),x}^*\omega\ ds\ dp\ .$$
With respect to the Riemannian metrics, the coordinate transformation
$$\SS^{n-1}\times\R^+_0\ni(p,s)\overset{y_x}{\longmapsto}y=\phi_{p,s}(x)\in M $$
has determinant
$$\det\frac{\partial y_x(p,s)}{\partial(p,s)}=2^{-m}\sinh(2s)^m\sinh(s)^{\dim M-m-1}\ .$$
Since $s=d(x,y_x(p,s))$, the transformation formula yields the following.
\begin{thm}\label{thmnegcurv}
Let $M$ be a negatively curved symmetric space and let $\omega\in\Omega^{k+1}$ be closed, i.e. $\d\omega=0$, and of compact support. Then we have $\omega=\d\mu$ with $\mu\in\Omega^k(M)$ defined by
\[{musymsp}\mu_x=\frac{-1}{\vol{\SS^{n-1}}}\intd_M \frac{2^m\ \L_{yx}^*\iota_{T_{yx}}\P_{yx}^*\omega}{\sinh(2\d(x,y))^m\sinh(\d(x,y))^{\dim M-m-1}}\ \ dy\ ,\]
and
\[{Lnegcurvpp}\L_{yx}=e^{-d(y,x)\sqrt{J_{yx}}}=\pmx{
1&0    &0\\
0&e^{-2d(y,x)}\id_m   &0\\
0&0&e^{-d(y,x)}\id_{\dim M-m-1}
}\]
as in \eqref{Lnegcurv}.
\end{thm}

\section{Compact Symmetric Spaces of Rank 1}\label{cross}
\subsection{Spheres}\label{sec:sn}
On the standard sphere $M=\SS^n$ of constant sectional curvature 1, for $p\in\SS^n$, we consider the local flow $\phi_{p}$ defined on
$$
\set{(x,t)}{x\neq\pm p, -d(x,-p)<t<d(x,p)}\subset\SS^n\times\R$$
so that for $v\in T_pM$, $\norm{v}=1$ we have 
$$\phi_{p,t}(\exp_p(sv))=\exp_p((s-t)v)\ .$$
We compute \eqref{flowmain} with $b(x)=d(x,p)$,
\[{muspherefst}\mu_{p,x}=-\intd_0^{d(x,p)}(\iota_{T_p}\phi_{p,t}^*\omega)_x\ dt\qmbox{for}d(x,p)<\pi\qmbox{i.e.}x\neq -p\ ,\]
and have
$$\d \mu_p=\omega$$
because
$$(\phi_{p,d(x,p)}^*\omega)(x)=\lim_{t\to d(x,p)}(\phi_{p,t}^*\omega)(x)=0\ .$$
As before we express the differential of the flow $\phi_p$ in terms of the parallel transport,
\spll{dfdfdfdfdf}{
d_x\phi_{p,t}T_{px}&=\P_{\phi_{p,t}(x),x}T_{px}\\
d_x\phi_{p,t}V&=\frac{\sin(\d(x,p)-t)}{\sin(\d(x,p))}\P_{\phi_{p,t}(x),x}V
}
for $V\in T_x\SS^n$, $V\perp T_{px}$, and $d(x,y)<\pi$, i.e. $x\neq -p$.
Since $\iota_{T_{p}}\omega$ annihilates $T_{p}$ in \eqref{muspherefst}, we only need the second equation in \eqref{dfdfdfdfdf} and get
\[{muspherescnd}\mu_{px}=-\intd_0^{\d(x,p)}\l(\frac{\sin(\d(x,p)-t)}{\sin(\d(x,p))}\r)^k  \iota_{T_{px}}\P_{\phi_{p,t}(x),x}^*\omega\ dt\qmbox{for}d(x,p)<\pi\ .\]
We remove the singularity at $x=-p$ by averaging over $p\in\SS^n$, 
$$\mu_x=\frac{-1}{\vol\SS^n}\intd_{\SS^n}\intd_0^{\d(x,p)}\l(\frac{\sin(\d(x,p)-t)}{\sin(\d(x,p))}\r)^k  \iota_{T_{px}}\P_{\phi_{p,t}(x),x}^*\omega\ dt\ d^np\ .$$
To bring this into the standard form \eqref{mfdinv}, we use the coordinate transformation
$$y=y(p,t)=\phi_{p,t}(x)\qmbox{,}s=d(x,p)\qmbox{,}d(x,y)=t\ .$$
With respect to the Riemannian metrics, this has the determinant
$$\det\frac{\partial(p,t)}{\partial(y,s)}=\l(\frac{\sin(s)}{\sin(d(x,y))}\r)^{n-1}\ .$$
The transformation formula yields
\spl{
\mu_x&=\frac{-1}{\vol\SS^n}\intd_{\SS^n}\intd_{d(x,y)}^{\pi}\l(\frac{\sin(s-\d(x,y))}{\sin(s)}\r)^k\l(\frac{\sin(s)}{\sin(d(x,y))}\r)^{n-1}\   \iota_{T_{yx}}\P_{yx}^*\omega\ ds\ d^ny\\
&=\frac{-1}{\vol\SS^n}\intd_{\SS^n}\frac{\intd_{d(x,y)}^{\pi}\sin(s)^{n-1-k}\sin(s-d(x,y))^k \ ds\ }{\sin(d(x,y))^{n-1}}\ \iota_{T_{yx}}\P_{yx}^*\omega\ d^ny\ .
}
\begin{thm}\label{dinvspherethm}
Let $\omega\in\Omega^{k+1}(\SS^n)$ be exact (equivalently closed and $k<n-1$ or $k=n-1$ and $\int_{\SS^n}\omega=0$). Let $\mu\in\Omega^k(\SS^n)$ be defined by
$$\mu_x=\intd_{\SS^n}\lambda(d(x,y))\  \iota_{T_{yx}}\P_{yx}^*\omega\ d^ny $$
where
\[{lambdasphere}\lambda(d)=\frac{-1}{\sin(d)^{n-1}\vol\SS^n}\intd_{d}^{\pi}\sin(s)^{n-1-k}\sin(s-d)^k \ ds\ .\]
\end{thm}
\proof We still need to show that $\d\mu=\omega$ if the form $\mu\in\Omega^k(M)$ is given by $\mu_x=\int_{\SS^n}\mu_{px}\ dp^n$ with $\mu_{px}$ as in \eqref{muspherescnd}. To that end we will show that $\int_C\mu=\int_D\omega$ for every embedded oriented $(k+1)$-disk $D\subset\SS^n$ with boundary $C=\partial D$. To show this, let $C_x$ denote the oriented volume element of $C$ at $x\in C$ and compute
\[{imx76}\intd_C\mu=\intd_C\l(\intd_{\SS^n}\mu_p\ dp\r)=\intd_C\intd_{\SS^n}\mu_{px}(C_x)\ dp\ dx=\intd_{\SS^n}\intd_C\mu_{px}(C_x)\ dx\ dp\ .\]
We may interchange the order of integration by Fubini's Theorem, since $\mu_{p,x}(C_x)$ is absolutely integrable on $C\times\SS^n$. This is because
$$\abs{\mu_{px}(C_x)}\leq\frac{K}{d(-p,x)^k}$$
for some constant $K$, and $k$ is less than the codimension of $C\cong\Delta C=\set{(-c,c)}{c\in C}=\set{(x,p)}{d(-p,x)=0}$ in $C\times\SS^n$.
\par
If $k<n-1$, then for almost all $p\in\SS^n$, we have $-p\not\in D$, hence $\mu_p$ is defined on $D$ by Cartan's formula and $\int_D\omega=\int_C\mu$ by Stokes' theorem. The same applies if $k=n-1$ and $p\not\in D$.
\par
If $k=n-1$ then still $p\not\in C$ for almost all $p$. If $p\in D$ then we close $C$ in $\SS^n\setminus{-p}\supset F$, $C=\partial F$, $F$ a $(k+1)$-disk in $\SS^n\setminus\st{-p}$. By assumption the integral of $\omega$ over $F\cup_C D$ vanishes. Since $F$ and $D$ induce opposite orientations on $C$ we have
$$\intd_D\omega=-\intd_F\omega=\intd_C\mu_p\ .$$
The second identity holds because $\d\mu_p=\omega$ on $\SS^n\setminus\st{-p}\supset F$, and Stokes' theorem.
\proofend

\subsection{Projective spaces $\CP^n$, $\HP^n$, $\OP^2$}\label{sec:cpn}
Let $M=\KP^n$, $\K=\C,\mathbb{H},\mathbb{O}$ be a projective space and $m=1,3,7$ respectively be the dimension of the imaginary part $\Im\K$ of $\K$.
The metric on these spaces will be normalised to diameter $\pi/2$, so that the sectional curvature of these spaces ranges from $1$ to $4$. Given $\omega\in\Omega^{k+1}(M)$ choose $r,s\in\N$ so that $r+s+1=n$ and
$$(m+1)r\leq k\leq (m+1)r+m\ .$$
Thus, in the case of $M=\CP^n$, we have $k=2r$ or $2r+1$. Let $\GM$ denote the projective Grassmannian, i.e.
the set of projective $r$-dimensional planes in $M$.  $\GM$ is naturally diffeomorphic with the usual Grassmanian $G_{r+1}(\K^{n+1})$ of linear $(r+1)$-dimensional subspaces of $\K^{n+1}$. The focal locus of $V\in\GM$ in $M$ is
$$\VV:=\set{x}{d(V,x)=\frac{\pi}{2}}\in G_s(\KP^n)\cong G_{s+1}(\K^{n+1}).$$
In the Grassmanian $G^\K_{r+1}(\K^{n+1})$ this corresponds to the orthogonal complement.
\medskip

For $V\in\GM$ we consider the local flow $\phi_V$ defined on
$$\set{(x,t)}{0<d(V,x)<\frac{\pi}{2},\ d(V,x)-\frac{\pi}{2}<t<d(V,x)}\subset M\times\R$$
so that for $p\in V$, $v\in N_p(V,M)$ a unit normal vector to $V$ and $0<s<\pi/2$ we have
\[{lkjlps}
\phi_{V,t}(\exp_psv)=\exp_p((s-t)v)\ .\]
The corresponding vector field is
$$T_V(x)=\diff{t}{0}\phi_{V,t}(x)\ .$$
As before, since $\lim_{t\to d(V,x)}(\phi_{V,t}^*\omega)(x)=0$ we get from \eqref{flowmain} that the form
\[{muposcurv}\mu_{V,x}=-\intd_0^{d(V,x)}\phi_{V,t}^*\iota_{T_V}\omega\ dt\qmbox{for}d(x,p)<\frac{\pi}{2}\qmbox{i.e.}x\not\in\VV\]
satisfies
$$\d\mu_V=\omega\ .$$
As in section \ref{secnegcurv} we compare $\phi_{V,s}$ to the parallel transport.
The orbits of the flow \eqref{lkjlps} are geodesics of length $\frac{\pi}{2}$ joining points $\pp\in\VV$ to $p\in V$.
We need the solutions of the Jacobi equation \eqref{diffeqjac} whose component parallel to $V$ vanishes at $\pp$, whose component parallel to $\VV$ vanishes at $p$, and whose component in the $\K$-span of $T_V$ vanishes at both $p$ and $\pp$.
The first two are eigenspaces of the Jacobi operator with eigenvalue $-1$.
The vector $T_V$ is an eigenvector with eigenvalue $0$ and the $\R$-orthogonal complement $\Im{\K}T_V$ of $T_V$ in $\K T_V$ is the eigenspace corresponding to the eigenvalue $-4$ of the Jacobi operator.
As in \eqref{Lnegcurv}, with $y=\phi_{V,t}(x)$, we get
\spll{Lposcurv}{\L_{yx}^V&:=\l(\P_{yx}\r)^{-1}d_x\phi_{V,d(y,x)}\\
&=\pmx{\frac{\cos(d(V,y))}{\cos(d(V,x))}\ \id_{2r}&0&0&0\\
0&1&0&0\\
0&0&\frac{\sin(2(d(V,y)))}{\sin(2d(V,x))}\ \id_m   &0\\
0&0&0&\frac{\sin(d(V,y))}{\sin(d(V,x))}\ \id_{2s}
}
}
the matrix being with respect to the splitting
$$T_xM=V\oplus\spn{T_V}\oplus\Im{\K}T_V\oplus\VV\ .$$
In order to obtain a form $\mu$ defined on all of $M$, we average \eqref{muposcurv} over $V\in\GM$,
\[{ujkd}\mu_x:=\frac{-1}{\vol\GM}\intd_{\GM} \intd_0^{d(V,x)}\phi_{V,t}^*\iota_{T_V}\omega\ dt\ dV\ .\]
To bring this into the form \eqref{mfdinv}, we use the coarea formula for the submersions
\[{coordtrafo}\arraycolsep2pt\def\arraystretch{1.1}
\begin{matrix}
&\GM\times\sg{0,\frac{\pi}{2}}&&\KP^n\times\sg{0,\frac{\pi}{2}}&\\
&\cup&&\cup&\\
F_x\colon&\set{(V,t)}{0<t<d(V,x)}&\longrightarrow&\set{(y,s)}{d(x,y)<s<\frac{\pi}{2}}&\\
&\upin&&\upin&\\
&(V,t)&\longmapsto&(\phi_{V,t}(x),d(V,x))&&
\end{matrix}
\]
The submersions $F_x$ are equivarient under isometries of $\KP^n$ fixing $x$. The fibres of $F_x$,
\spl{
F_x^{-1}(y,s)&=\set{(V,t)}{\phi_{V,t}(x)=y,\ d(V,x)=s}\\
&=\set{V}{\gamma_{yx}(s)\in V\perp_\K \gamma_{xy}'(s)}\times\st{d(x,y)}
}
are extrinsically homogeneous, the group of isometries of $\GM$ fixing $x$ and $y$ acts transitively on $F_x^{-1}(y,s)$. Therefore the Jacobian
$$\Jac:=\sqrt{\det (d_{V,t} F_x) (d_{V,t} F_x)^*}$$
for any $(V,t)$ with $F_x(V,t)=(y,s)$ is well defined as a function of $(y,s)$. By the coarea formula (see \cite{fd}, Thm. 3.2.12 for instance) for the submersions $F_x$, we can rewrite the average \eqref{ujkd} as
\spl{
\mu_x&=\frac{-1}{\vol\GM}\intd_{\KP^n} \intd_{d(x,y)}^{\frac{\pi}{2}}\frac{\iota_{T_{yx}}}{\Jac}\intd_{F_x^{-1}(y,s)}\Lambda^kd_x\phi_{V,t}^*\ dt\ dV\ ds\ \omega\ dy\\
&=\frac{-1}{\vol\GM}\intd_{\KP^n} \intd_{d(x,y)}^{\frac{\pi}{2}}\frac{\iota_{T_{yx}}}{\Jac}\intd_{\overset{V\ni\gamma_{yx}(s)}{V\perp_\K \gamma_{xy}'(s)}}\Lambda^kd_x\phi_{V,d(x,y)}^*\ dV\ ds\ \omega\ dy\ .
}

Since $d_x\phi_{V,d(x,y)}=\P_{yx} \L_{yx}^V$, this yields the following.
\begin{thm}
Let $\omega\in\d\Omega^k(\KP^n)\subset\Omega^{k+1}(\KP^n)$. Then $\omega=\d\mu$ with $\mu\in\Omega^k(\KP^n)$ given by
\[{maindinvsymposcurv}\mu_x
=\frac{-1}{\vol\GM}\intd_{\KP^n} \intd_{d(x,y)}^{\frac{\pi}{2}}\fracd{\intd_{\overset{V\ni\gamma_{yx}(s)}{V\perp_\K \gamma_{xy}'(s)}}\Lambda^k {\L_{yx}^V}^*\ dV}{\Jac}\ ds\ \iota_{T_{yx}}\P_{yx}^*\omega\ dy \]
where $\L_{yx}^V$ is given by \eqref{Lposcurv}.
\end{thm}

\proof As for the proof of Theorem \ref{dinvspherethm}, we still need to show that $\d\mu=\omega$ if the form $\mu\in\Omega^k(M)$ is given by $\mu=\int_{\GM}\mu_{V}\ dV$ with $\mu_{V}$ as in \eqref{muposcurv}, defined on $M\setminus\VV$. To that end we will show that $\int_C\mu=\int_D\omega$ for every embedded oriented $(k+1)$-disk $D\subset M$ with boundary $C=\partial D$. To show this, let $C_x$ denote the oriented volume element of $C$ at $x\in C$ and compute
\spll{imxcp}{
\intd_C\mu&=\intd_C\intd_{\GM}\mu_V\ dV=\intd_C\intd_{\GM}\mu_{V,x}(C_x)\ dV\ dx\\
&=\intd_{\GM}\intd_C\mu_{V,x}(C_x)\ dx\ dV\ .
}
In \eqref{imxcp} we may interchange the order of integration by Fubini's Theorem, since $\mu_{V,x}(C_x)$ is absolutely integrable over $\GM\times C$. This is because
$$\abs{\mu_{V,x}(C_x)}\leq\frac{K}{\l(\frac{\pi}{2}-d(V,x)\r)^{\dim_\R V +m}}=\frac{K}{\l(d(\VV,x)\r)^{\dim_\R V +m}}$$
for some constant $K$. This is integrable over a tubular neighbourhood of $\VV$ since the codimension of $\VV$ is $\dim_\R V +m+1$. 
\par
The map
$$\arraycolsep2pt\def\arraystretch{1.1}
\begin{matrix}
&\set{(V,v)}{V\in\GM,v\in SN_p(V,\KP^n), p\in V}&\longrightarrow&\KP^n&\\
&\upin&&\upin&\\
&(V,v)&\longmapsto&\exp_p\l(\frac{\pi}{2}v\r)&&
\end{matrix}
$$
is a submersion, and maps $SN(V,\KP^n)\to\VV$. By parametric transversality, for almost all $V\in\GM$, $\VV$ is therefore transversal to $D$ (and hence $C$).
This means that for almost all $V\in\GM$, the intersection $\VV\cap C=\emptyset$, and $\VV\cap D$ is empty, or possibly finite in case $k=\dim_\R V+m$.
\par
In the first case, we have $\omega=\d\mu_V$ on all of $D$ because $D$ lies in the domain of definition of $\mu_V$. 
In the second case, we can close $C=\partial W$ with a domain $W$ not intersecting $\VV$. Thus $\mu_V$ is defined on $W$ and again by Stokes' theorem, observing that $W$ and $D$ induce opposite orientations on $C$, and $\int_{W\cup_{C}D}\omega=0$, we get
$$\intd_C\mu_V=-\intd_W\omega=\intd_D\omega\ .$$ 
\proofend

We compute an example to illustrate the application of \eqref{maindinvsymposcurv}. Let $\omega\in\d\Omega^1(\CP^2)\subset\Omega^2(\CP^2)$. We then have $k=1$, $r=0$, $s=1$. We need the operator $\L_{yx}^V$ of \eqref{Lposcurv} only on $T_{yx}^\perp$, where, with respect to the splitting
\[{cp2splt1}T_{yx}^\perp=\spn{iT_{yx}}\oplus T^{\perp_\C}\subset T_x\CP^2\ ,\]
it is given by the $(3\times 3)$-matrix 
$$\L_{yx}^V=\pmx{
\frac{\sin(2(d(V,y)))}{\sin(2d(V,x))}   &0\\
0&\frac{\sin(d(V,y))}{\sin(d(V,x))}\ \id_{2}
}\ .
$$
The innermost integral in \eqref{maindinvsymposcurv} is over one point only and a straightforward calculation gives
\spl{
\frac{1}{\Jac}&=\abs{\det\pmx{
0&1&0&0&0\\
1&0&0&0&0\\
0&0&\frac{\sin(2s)}{\sin(2d(x,y))}&0&0\\
0&0&0&\frac{\sin(s)}{\sin(d(x,y))}&0\\
0&0&0&0&\frac{\sin(s)}{\sin(d(x,y))}
}}\\
&=\frac{\sin(2s)\sin^2(s)}{\sin(2d(x,y))\sin^2(d(x,y))}\ .
}
Inserting this into \eqref{maindinvsymposcurv} we obtain
\spll{cp2dinv}{\mu_x&=\frac{-1}{\vol\CP^2}\intd_{\CP^2} \intd_{d(x,y)}^{\frac{\pi}{2}}\frac{\iota_{T_{yx}}}{\Jac}\l(\int_{p_{yx}\in V\perp_\K T}\Lambda^k \L_{yx}^V\ dV\r)^*\P_{yx}^*\omega\ ds\ dy\\
&=\frac{-2}{\pi^2}\intd_{\CP^2} \frac{\iota_{T_{yx}}}{\sin(2d(x,y))\sin^2(d(x,y))} \L_{yx}^*\P_{yx}^*\omega\ dy\ .
}
The matrix of $\L_{yx}$ with respect to the splitting \eqref{cp2splt1} is
\spll{cptwocoefs}{
\L_{yx}&=\intd_{d(x,y)}^{\frac{\pi}{2}}\sin(2s)\sin^2(s)\pmx{
\frac{\sin(2(s-d(x,y)))}{\sin(2s))}   &0\\
0&\frac{\sin(s-d(x,y))}{\sin(s)}\ \id_{2}
} \ ds\\
&=\intd_{d(x,y)}^{\frac{\pi}{2}}\pmx{
\sin(2(s-d(x,y))\sin^2(s)&0\\
0&\sin(s-d(x,y))\sin(2s)\sin(s)\ \id_{2}
} \ ds\\
&=\pmx{L^0_{yx}&0\\
0&L^1_{yx}\ \id_{2}
}
}
with
\spl{
L^0_{yx}&=\frac{\l(\pi-2d(x,y)\r)\,\sin(2d(x,y))+2\cos(2d(x,y))+2}{8}\\
L^1_{yx}&=\frac{\cos(3d(x,y))+(4d(x,y)-2\pi)\sin(d(x,y))+7\cos(d(x,y))
 }{16}
}

\section{Linking Numbers}\label{sec:link}
In terms of differential forms, the linking number of two closed oriented submanifolds $K,L\subset M$ of an oriented manifold $M$, whose fundamental classes are trivial in $M$, can be defined as
$$\lk{K,L}=\intd_Kd^{-1}\eta_L $$
where $\eta_L$ is a Poincar\'e dual of $L$, i.e. any $\eta_L\in\Omega^{m-l}(M)$, so that
$$d\eta_L=0\qmbox{and}\intd_L \alpha=\intd_M \alpha\wedge\eta_L$$
for all $\alpha\in\Omega^l(M)$, $d\alpha=0$, see \cite{botu}.
The Poincar\'e dual can be chosen to have support in any tubular neighbourhood of $L$.
Thus every integral formula for an inverse of the Cartan derivative gives a Gauss-type integral for the linking number.
\begin{thm}\label{thm:link}
Let $M$ be a connected oriented Riemannian manifold and $\L\in\Hom(TM,TM)$ and $\lambda\in C^\infty(\R^+)$ are so that $\eqref{mfdinv}$ holds.
Let $K,L\subset M$ be connected closed oriented submanifolds of dimensions $k,l$ respectively, with $k+l+1=\dim M$.
If the fundamental classes of $K$ and $L$ vanish in $M$, then the linking number of $K$ with $L$ is
\[{lkmain}\lk{K,L}=\intd_K\intd_L \lambda(d(x,y))\vol^M_y(\d L_y\wedge\P_{yx}(T_{yx}\wedge\Lambda^k\L_{yx}\d K_x))\ d^ly\  d^kx\ .\]
\end{thm}
\proof We have \eqref{mfdinv} for the Poincar\'e dual $\eta_L$ of $L$,
$$\eta_L=d\mu$$
$$\mu_x=\intd_M \lambda(d(x,y)) \L_{yx}^*\iota_{T_{yx}}\P_{yx}^*\eta_L(y)\ d\vol_g^M(y)\ .$$
Let $\d L_y$, $\d K_x$ denote the volume elements of $L$ and $K$ at $y\in L$ respectively $x\in K$.
Since we can choose the support of $\eta_L$ in an arbitrary small tubular neighbourhood of $L$, we can compress the integral over $M$ to one over $L$, as in \cite{SS}. Thus
\spl{
\lk{K,L}&=\intd_K\mu_x(\d K_x) d^kx\\
&=\intd_K\intd_M \lambda(d(x,y)) \L_{yx}^*\iota_{T_{yx}}\P_{yx}^*\eta_{L,y}(\d K_x)\ d^ny\  d^kx\\
&=\intd_K \intd_L \lambda(d(x,y)) \L_{yx}^*\iota_{T_{y,x}}\P_{yx}^*\eta_{L,y}(\d K_x)\ d^ly\  d^kx\\
&=\intd_K\intd_L \lambda(d(x,y))\eta_{L,y}(\P_{yx}(T_{yx}\wedge\Lambda^k\L_{yx}\d K_x))\ d^ly\  d^kx\\
&=\intd_K\intd_L \lambda(d(x,y))\vol^M_y(\d L_y\wedge\P_{yx}(T_{yx}\wedge\Lambda^k\L_{yx}\d K_x))\ d^ly\  d^kx
}
\proofend

\subsection{Example: negatively curved symmetric spaces}
From Theorems \ref{thmnegcurv} and \label{thm:link} we immediately get a formula for the linking number in the negatively curved symmetric spaces $M=H^n,\CH^n,\HH^n,\OH^2$ with $m=0,1,3,7$ respectively. If $K^k,L^l\subset M$ are oriented closed submanifolds with $k+l+1=(m+1)n=\dim M$, then the linking number is
\[{lknegcursym}\lk{K,L}=\frac{-1}{\vol{\SS^{n-1}}}\intd_K\intd_L\frac{2^m\ \vol_y^M(\d L_y\wedge\P_{yx}(T_{yx}\wedge\Lambda^k \L_{yx}\d K_x))}{\sinh(2\d(x,y))^m\sinh(\d(x,y))^{\dim M-m-1}} \ d^ly d^kx\ .\]
A matrix for $\L_{yx}$ is explicitly given in \eqref{Lnegcurvpp}. In the case of the hyperbolic space $M=H^n$, we have $m=0$ and $\L_{yx}$ is multiplication with $e^{-d(x,y)}$. The linking integral \eqref{lknegcursym} simplifies to
$$\lk{K,L}=\frac{-1}{\vol{\SS^{n-1}}}\intd_K\intd_L\frac{e^{-kd(y,x)}\vol_y^M(\d L_y\wedge\P_{yx}(T_{yx}\wedge\d K_x))}{\sinh(\d(x,y))^{\dim M-m-1}} \ d^ly d^kx\ . $$

\subsection{Example: $M=\SS^n$}\label{exlksn}
In this case \eqref{lkmain} becomes
$$\lk{K,L}=\intd_K\intd_L\lambda(d(x,y)) \vol_y^M(\d L_y\wedge\P_{yx}(T_{yx}\wedge\d K_x))\ d^ly d^kx\ .$$
where $\lambda$ is as defined in \eqref{lambdasphere}.
We orient the sphere by the outward normal vector, i.e. so that for $y\in\SS^n\subset\R^{n+1}$ we have
$$\vol_y^M(w_1\ldots,w_n)=\det(w_1,\ldots,w_n,y)$$
for $w_1,\ldots,w_n\in T_y\SS^n=y^\perp\subset\R^{n+1}$.
Then
\spl{\vol_y^M(\d L_y\wedge\P_{y,x}(T_{y,x}\wedge\d K_x))&=\frac{\det(L_1,\ldots,L_l,-x,K_1,\ldots,K_k,y)}{\sin(d(x,y))}\\
&=\frac{[dL,x,dK,y]}{\sin(d(x,y))}
}
in the notation of \cite{dtg3}. Thus
$$\lk{K,L}=\intd_K\intd_L\frac{\lambda(d(x,y))}{\sin(d(x,y))}\det(\d L_y,-x,\d K_x,y)\ dxdy$$
as computed in \cite{dtg3}.

\subsection{Example: $M=\CP^2$}
Let $K^1,L^2\subset\CP^2$ be connected closed oriented nullhomologous submanifolds of dimension $1$ and $2$ respectively. From \eqref{lkmain} and \eqref{cp2dinv},
$$\lk{K,L}=\frac{-2}{\pi^2}\intd_K\intd_L\frac{\vol^{\CP^2}_y(\d L_y\wedge\P_{y,x}(T_{y,x}\wedge\L_{y,x}\d K_x))}{\sin(2d(x,y))\sin^2(d(x,y))} \ d^ly\  d^kx\ .$$
Splitting $\d K_x=\kappa_0(x)iT+\kappa_1$ and $\d L_y=\lambda_0\wedge iT+\lambda_1$ into components perpendicular respectively parallel to $iT_{yx}$, this becomes
\spl{
\lk{K,L}=\frac{-2}{\pi^2}\intd_K\intd_L&\L^0_{y,x} \vol^{\CP^2}_y(\lambda_1(y)\wedge\P_{y,x}(T_{y,x}\wedge\kappa_0(x)iT))\\
&+\L^1_{y,x} \vol^{\CP^2}_y(\lambda_0(y)\wedge iT\wedge\P_{y,x}(T_{y,x}\wedge\kappa_1(x)))\ d^ly\  d^kx\ .
}
where, because of \eqref{cptwocoefs}, the coefficients are
\spl{
\L^0_{y,x}&=-\frac{\l(2\,d(x,y)-\pi\r)\,\sin(2d(x,y))-2\cos(2d(x,y))-2}{8\sin(2d(x,y))\sin^2(d(x,y))}\\
\L^1_{y,x}&=\frac{\cos(3d(x,y))+(4d(x,y)-2\pi)\sin(d(x,y))+7\cos(d(x,y))
 }{16\sin(2d(x,y))\sin^2(d(x,y))}\ .
}

\affiliation{Department of Mathematics and Statistics, Maynooth University, Maynooth, Co. Kildare, Ireland}\\
\email{stefan.bechtluft-sachs@mu.ie}
\medskip

\affiliation{University of Cyprus, Department of Mathematics and Statistics, P.O.Box 20537, 1678 Nicosia, Cyprus}\\
\email{samiou@ucy.ac.cy}

\end{document}